\documentclass[12pt, reqno]{amsart}

\usepackage{latexsym, amsfonts, amsmath, amsthm, amscd} 
\newcommand{\pair}[2]{(#1|#2)}

\newcommand{\apair}[2]{\langle#1|#2\rangle}
\newcommand{\apaircd}{\langle\cdot|\cdot\rangle}

\newcommand{\tsum}{\textstyle\sum}

\newcommand{\boxt}{\Tox\kern -6.3pt\raise .55pt
    \hbox{$\scriptstyle{\times}$}}

\newcommand{\bd}{{\mathbf d}}

\newcommand{\bh}{{\mathbf h}}

\newcommand{\bp}{{\mathbf p}} 
\newcommand{\bq}{{\mathbf q}}  
  
\newcommand{\bfr}{{\mathbf r}}

\newcommand{\D}{{\mathfrak D}}

\newcommand{\ltimes}{\vbox to 5.4pt{\leaders\vrule\vfil}\kern
  -5pt\times}

\newcommand{\del}{\partial}

\newcommand{\llongrightarrow}{\relbar\joinrel\longrightarrow}

\newcommand{\mmapright}[1]{\;\smash{\mathop
   {\llongrightarrow}\limits^{#1}}\;}

\newtheorem{thm}{Theorem}[section]
\newtheorem{lem}[thm]{Lemma}

\newtheorem{prop}[thm]{Proposition}
\newtheorem{cor}[thm]{Corollary}

\theoremstyle{definition}

\newtheorem{remark}[thm]{Remark}

\numberwithin{equation}{section}

\newcommand{\inv}{^{-1}} 
\newcommand{\ovl}[1][]{\overline{#1}}

\newcommand{\half}[1][1]{\frac{#1}{2}}
\newcommand{\thalf}{\textstyle\frac{1}{2}}

          \newcommand{\g}{{\mathfrak g}}

\newcommand{\fsu}{\mathfrak{su}} 

\newcommand{\fsl}{\mathfrak{sl}}

\newcommand{\ophi}{\ovl[\phi]}
\newcommand{\opsi}{\ovl[\psi]}

\newcommand{\ox}{\ovl[x]} 

\newcommand{\oA}{\ovl[A]} \newcommand{\oB}{\ovl[B]}

\newcommand{\oa}{\ovl[a]}

\newcommand{\C}{{\mathbb C}}

\newcommand{\R}{{\mathbb R}}
  
 \newcommand{\bbT}{{\mathbb T}} 

\newcommand{\gR}{\g_{\R}}
\newcommand{\GR}{G_{\R}}  \newcommand{\XR}{X_{\R}}
\newcommand{\Gc}{G_{c}}  \newcommand{\gc}{\g_{c}}

  \newcommand{\PR}{P_{\R}} 
 
\newcommand{\OR}{\OO_{\R}}
\newcommand{\muR}{\mu_{\R}}

\newcommand{\al}{\alpha}
\newcommand{\be}{\beta}
\newcommand{\ga}{\gamma}

\newcommand{\la}{\lambda} \newcommand{\La}{\Lambda}

\newcommand{\sig}{\sigma}

                                     \newcommand{\A}{{\mathcal A}}
                                     \newcommand{\B}{{\mathcal B}} 
 
\newcommand{\cD}{{\mathcal D}} 
\newcommand{\E}{{\mathcal E}} 
                                          \newcommand{\F}{{\mathcal F}}
\newcommand{\cG}{{\mathcal G}}
\newcommand{\cH}{{\mathcal H}}

\newcommand{\cL}{{\mathcal L}}

                                          \newcommand{\OO}{{\mathcal O}} 
\newcommand{\cP}{{\mathcal P}} 
                                            \newcommand{\Q}{{\mathcal Q}}

                                         \newcommand{\T}{{\mathcal T}}  
                                            
\newcommand{\cV}{{\mathcal V}}

\newcommand{\Pol}{\mathop{\mathrm{Pol}}\nolimits}

\newcommand{\gr}{\mathop{\mathrm{gr}}\nolimits}

\newcommand{\End}{\mathop{\mathrm{End}}\nolimits}

\newcommand{\gog}{\g\oplus\g}

\newcommand{\Ug}[1][]{{\mathcal U}_{#1}({\mathfrak g})} 
 
\newcommand{\Sg}[1][]{S^{#1}({\mathfrak g})}

\newcommand{\CP}{\mathbb{CP}}
\newcommand{\RP}{\mathbb{RP}}

\newcommand{\ssig}{^{\sig}}

\newcommand{\gs}{\g^{\sharp}}

\newcommand{\TX}{T^*X} 
\newcommand{\TXR}{T^*\XR} 
\newcommand{\RTX}[1][]{R^{#1}(T^*X)}
 
\newcommand{\PolTXR}[1][]{\Pol^{#1}(T^*\XR)}

\newcommand{\smhalf}{{\scriptscriptstyle \frac{1}{2}}}
\newcommand{\DXR}{\D^{\smhalf}(\XR)}
\newcommand{\DXsmooth}{\D^{\smhalf}(X)}

 \newcommand{\DX}[1][\smhalf]{\D^{#1}_{alg}(X)}

\newcommand{\dhalf}{\delta^{\smhalf}}

\newcommand{\Ehalf}[1][\smhalf]{\E^{#1}} 
 
\newcommand{\etah}{\eta_{\smhalf}}
\newcommand{\xih}{\xi_{\smhalf}}
 
\newcommand{\GEhalf}{\Gamma(X,\Ehalf_X)}
\newcommand{\GEhalffin}{\GEhalf^{\Gc-fin}}

\newcommand{\sal}{^{\al}}
\newcommand{\sbe}{^{\be}}
 
\newcommand{\rep}{representation}  
\newcommand{\reps}{representations }

\newcommand{\CR}{C^{\scriptscriptstyle\RR}}

\newcommand{\RR}{{\mathcal R}}  \newcommand{\RRm}{{\mathcal R}_{\mu}}
   \newcommand{\DD}{{\mathcal D}}

\newcommand{\hatRR}{\widehat{\RR}}

\newcommand{\shol}{^{hol}}
\newcommand{\Ahol}{A\shol} 

\newcommand{\Afin}{\A_{fin}} \newcommand{\Bfin}{\B_{fin}}

\newcommand{\Polhol}[1][]{\Pol_{hol}^{#1}(T^*X) }

\begin{document}   

\title 
{Equivariant Deformation Quantization for  the   Cotangent Bundle  of  a
Flag Manifold}
\author{Ranee Brylinski}
\address{Department of Mathematics,
         Penn State University, University Park 16802}
\email{rkb@math.psu.edu}
\urladdr{www.math.psu.edu/rkb}

\begin{abstract}
Let $\XR$ be a (generalized) flag manifold of a non-compact real semisimple Lie 
group $\GR$, where   $\XR$ and $\GR$ have complexifications $X$ and $G$. We
investigate the problem of constructing a graded   star product   on  
$\Pol(T^*\XR)$ which corresponds to a $\GR$-equivariant  quantization of symbols
into   smooth differential operators acting  on half-densities on $\XR$. 

We show that any  solution is algebraic in  
that it  restricts to  a $G$-equivariant graded   star product $\star$ on
the algebraic part $\RR$ of  $\Pol(T^*\XR)$. 
We construct, when  $\RR$ is
generated by the momentum functions  $\mu^x$ for $G$,  a preferred
choice of  $\star$ where $\mu^x\star\phi$ has the  form
$\mu^x\phi+\half\{\mu^x,\phi\}t+\Lambda^x(\phi)t^2$. 
Here $\Lambda^x$ are  
operators on $\RR$ which are  not  differential in the known examples and  
so   $\mu^x\star\phi$ is not local in $\phi$.

$\RR$ acquires  an invariant   positive definite  inner product
compatible with its  grading. The completion of  $\RR$ is a new  
Fock space type  model  of the   unitary representation of $G$ on $L^2$
half-densities on $X$.
\end{abstract}

\maketitle

\section{Introduction} 
\label{sec_intro}      
The equivariant deformation quantization (EDQ) problem for cotangent bundles is to
construct a graded $\GR$-equivariant star product $\star$ on the symbol algebra
$\PolTXR$ where $\XR$ is a homogeneous space of a real Lie group $\GR$.
We  require  that specialization of $\star$ at $t=1$ 
produces  the  algebra $\D(\XR,\cL)$ of   smooth differential operators 
for some $\GR$-homogeneous line bundle $\cL$. 
Then  $\star$ corresponds to a   quantization map 
$\Q$ from $\PolTXR$ onto $\D(\XR,\cL)$; $\GR$-equivariance of $\star$
amounts to  $\GR$-equivariance  of $\Q$.

Motivated by geometric
quantization (GQ), we   take $\cL$ to be the half-density line bundle
$\Ehalf_{\XR}$.  Let $\DXR=\D(\XR,\Ehalf_{\XR})$.  This choice of $\cL$ is 
naturally consistent with our requiring parity for $\star$.

If $\GR$ is compact,  the geometric methods of Fedosov  should admit
an equivariant version which leads to a
positive solution to this problem;  see  e.g., \cite{B-N-P-W}.
The resulting star product would be local, i.e. bidifferential, and so 
would extend to the full algebra of smooth functions on $\TXR$.

If $\GR$ is not compact,   the situation is very different. The known geometric
methods break down and we expect there is no bidifferential solution in general.
For instance, when $\GR=SL_{n+1}(\R)$ and $\XR=\RP^n$ ($n\ge 1$),
Lecomte and Ovsienko constructed  in \cite{L-O}     
a unique solution for $\Q$.  
The corresponding star product on $\Pol(T^*\RP^n)$ is not local, but locality is
violated in a nicely controlled way; see \cite{me:RPn}.  
We  imagine this  star product  will ultimately be ``explained" by
some new  non-local quantization scheme set in a more general framework.

In this paper, we investigate the EDQ problem for $\TXR$ when $\GR$ is a
non-compact real semisimple Lie group  and $\XR$ is a  flag manifold of
$\GR$. We assume that $\GR$ and $\XR$ have 
complexifications $G$ and $X$.  
Flag manifolds are the most familiar compact homogeneous spaces of $\GR$;
they exemplify the phenomenon of a  big symmetry group acting on a small  space.
The existence of $\Q$ is known 
in only the two cases: when $\XR$ is $\RP^n$ as discussed above or
(\cite{D-L-O1})   when $\XR$  is  the projectivised cone of null vectors in 
$\R^{p,q}$, $p+q\ge 5$, and $\GR=SO(p,q)$. 

We show in \S\ref{sec_RR} that any solution $Q$ is  algebraic in the sense that
it maps the algebraic part $\RR$ of $\PolTXR$ onto the
algebraic part $\DD$ of $\DXR$.
In the ``good case" (for instance when $G=SL_n(\C)$ or if $X$ is the full flag
manifold)  $\RR$ is generated by the momentum functions  $\mu^x$ 
where  $x$ lies in $\g=Lie(G)$. Then  $\RR=\Sg/I$   and  $\DD=\Ug/J$.

So by restriction  any solution $\Q$ defines a  $G$-equivariant quantization map
$\bq:\RR\to\DD$.  It is easy  to describe all such maps  $\bq$ (\S\ref{sec_exist}).
This suggests that we  construct a \emph{preferred} choice of
$\bq$, i.e., one that is special in some way, and then try to extend $\bq$ to $\Q$,
or equivalently, extend the corresponding star product on $\RR$ to $\PolTXR$.
See \S\ref{sec_RR} and  \S\ref{sec_Lax} for some preliminary ideas on the extension
problem.

We construct a preferred choice of $\bq$ in Theorem \ref{thm:main}, for the good
case.  Here is our method. Results in \rep\ theory 
of Conze-Berline and Duflo (\cite{CB-D}) and   Vogan (\cite{Vo1984})
give     a canonical  embedding   
$\Delta$ of  $\DD$ into the space  of smooth half-densities on $X$
(\S\ref{sec_I}); here we regard $X$ as a real manifold.  
We give a new geometric formula for $\Delta$ in (\ref{eq:Delta}).
The natural pairing $\int_X\al\ovl[\be]$ of half-densities
induces a  positive definite  inner product $\ga$ on $\DD$.
The $\ga$-orthogonal splitting of the order filtration on $\DD$ defines
our   $\bq$. 

In this way, $\RR$ acquires a positive definite inner product
$\apair{\phi}{\psi}=\ga(\bq(\phi),\bq(\psi))$ where the grading of $\RR$
is  orthogonal. 
Then $\apaircd$ is new even if  $\bq$ was unique to begin with (so if the 
\rep\ of $G$ on $\RR$ is multiplicity free).
The completion of  $\RR$ is a new  
Fock space type  model  of the   unitary representation of $G$ on $L^2$
half-densities on $X$ (\S\ref{sec_fock}).

Now $\bq$ defines a preferred graded $G$-equivariant star product $\star$ on
$\RR$. We find  in Corollary \ref{cor:3term}
that the star product $\mu^x\star\phi$ of a momentum function
with an arbitrary function in $\RR$ has the form 
$\mu^x\phi+\half\{\mu^x,\phi\}t+\La^x(\phi)t^2$ where $\La^x$ is the 
$\apaircd$-adjoint of ordinary multiplication by $\mu^{\sig(x)}$ ($\sig$ is a
Cartan involution of $\g$). 
This property that $\mu^x\star\phi$ is a three term sum  uniquely determines  
$\bq$ (Proposition \ref{prop:sat}).  
The $\La^x$ completely determine $\star$, but
they    are   not  differential in the known examples; see 
\S\ref{sec_Lax}.   Thus  $\mu^x\star\phi$ is not local in $\phi$.
  
An important feature is that $\DD$ has a natural trace functional $\T$
(Proposition \ref{prop:T}).
We give a  formula computing $\T$ by integration in (\ref{eq:T}). Then 
$\apair{\phi}{\psi}=\T(\bq(\phi)\bq(\psi\ssig))$ where $\sig$ is some anti-linear
involution of $\RR$; see (\ref{eq:Tcirc}).  

The philosophy here  is that the irreducible unitary \rep\ of $G$ on $L^2$
half-densities on $X$, modeled on $\RR$, should occur at the root of 
a solution to the  EDQ problem for $\TXR$. (This is certainly true when $\bq$ is
unique.)      The properties of this model, 
in particular the  interaction between   the Poisson algebra structure on
$\RR$ and  the inner product $\apaircd$,  should  control if and how $\star$ 
extends from $\RR$ to $\PolTXR$.  
 
I thank Pierre Bieliavsky, Jean-Luc Brylinski, Michel Duflo, Christian Duval, 
Simone Gutt, Valentin Ovsienko, Stefan Waldmann, and Alan Weinstein
for useful conversations last summer. I especially thank David Vogan for
discussions  in November 1999 which led to this paper.

\section{Cotangent bundles of flag manifolds}  
\label{sec_cot}   
Let $\GR$ be a non-compact connected real form
of a   complex semisimple Lie group $G$.
Let $\XR$   be a  (generalized) flag manifold   of $\GR$; then its 
complexification  $X$ is  a (generalized) flag manifold of $G$.
So $\XR=\GR/\PR$  and  $X=G/P$ are compact homogeneous spaces.
The classification of flag manifolds is well known. 

For example,  if  $\GR=SL_n(\R)$  then $G=SL_n(\C)$ and
their flag manifolds are $X^\bd(\R)$ and $X^\bd(\C)$ where 
$\bd=(d_1,\dots,d_s)$ with $1\le d_1<\cdots< d_s\le n-1$. Here
$X^\bd(\mathbb F)$  parameterizes the   flags   
$V=(V_1\subset\cdots\subset V_s)$  in
$\mathbb F^n$ where  $\dim V_j=d_j$. The simplest cases are the grassmannians 
of  $k$-dimensional subspaces in $\R^n$ and $\C^n$.
Here the flag 
manifolds of $SL_n(\R)$ and   $SL_n(\C)$ are in natural bijection; this happens
whenever  $\GR$ is the so-called split real form of $G$. 

The smooth action  of $\GR$ on $\XR$  lifts  canonically to a Hamiltonian action
on  $\TXR$ with moment map $\muR:\TXR\to\gR^*$.
Similarly, the holomorphic action  of $G$ on $X$  lifts  canonically to a 
Hamiltonian action on $\TX$ with  moment map
$\mu:T^*X\to\g^*$. Then $\muR=\mu|_{\TXR}$.
These moment maps embed  the cotangent spaces of $\XR$ and $X$  
into $\gR^*$ and $\g^*$. In our example, the cotangent space 
of   $X^\bd(\mathbb F)$ at $V$   identifies with the
subspace of $\fsl_n(\mathbb F)$  consisting of  maps 
$e:\mathbb F^n\to\mathbb F^n$ such that
$e(V_j)\subseteq V_{j-1}$.

Let $\PolTXR$ be the algebra of complex-valued smooth functions on
$T^*\XR$ which are polynomial on the cotangent fibers.
Then we have the algebra grading
\begin{equation}\label{eq:PolT*XR} 
\PolTXR=\oplus_{d=0}^\infty\PolTXR[d]
\end{equation}
by homogeneous degree along the fibers.
Clearly $\PolTXR$ is a graded Poisson algebra where
$\{\phi,\psi\}$ is homogeneous of degree $j+k-1$ if $\phi$ and $\psi$ are 
homogeneous of degrees $j$ and $k$.
We define $\phi\mapsto\phi\sal$ by $\phi\sal=(-1)^d\phi$ if  $\phi$ is
homogeneous of degree $d$; then  $\{\phi,\psi\}\sal=-\{\phi\sal,\psi\sal\}$.  
We define  $\phi\mapsto\ophi$ by  pointwise complex conjugation.

The Hamiltonian action of $\GR$ on $\TXR$ defines
a natural (complex linear) \rep\ of $\GR$   on $\PolTXR$. Then $\GR$ acts 
by graded Poisson algebra automorphisms which commute with $\al$ and
complex conjugation.
The  corresponding  \rep\ of $\g$    on $\PolTXR$ is given by the operators 
$\{\mu^x,\cdot\}$, $x\in\g$, where  $\mu^x\in\PolTXR[1]$ are the  
momentum functions. 

$\PolTXR$ is interesting  because  it is the algebra of symbols  
for (linear) differential operators acting on sections of a line bundle over $\XR$. 
 
\section{Equivariant star product problem for  $\TXR$}  
\label{sec_eqprob}

Our   motivating   problem  is  to construct a graded
$\GR$-equivariant star product  (with parity)  on  $\A=\PolTXR$.
This means that we want  an associative product $\star$ on $\A[t]$ which
makes $\A[t]$ into   an  algebra over $\C[t]$ in the following way.
If  $\phi,\psi\in\A$, then the product has the form
\begin{equation}\label{eq:star=} 
\phi\star\psi=\phi\psi+\thalf\{\phi,\psi\}t+
\tsum_{p=2}^\infty C_p(\phi,\psi)t^p
\end{equation}
where the  coefficients  $C_p$ satisfy
\begin{equation}\label{eq:star_list} 
\begin{array}{ll} 
{\rm (i)}  &  \mbox{$C_p(\phi,\psi)\in\A^{j+k-p}$ if
$\phi\in\A^j$  and $\psi\in\A^k$} \\[2pt] 
{\rm (ii)} &  C_p(\phi,\psi)=(-1)^p C_p(\psi,\phi)  \\[2pt]   
{\rm (iii)} & \ophi\star\opsi=\ovl[\phi\star\psi]  \\[2pt]
{\rm (iv)} & C_p(\cdot,\cdot)  \mbox{ is local on } \XR  \\[2pt]  
{\rm (v)} & \mu^x\star\phi-\phi\star\mu^x=t\{\mu^x,\phi\} \mbox{ for all $x\in\g$}
\end{array}
\end{equation}
\noindent 
 
Axiom (ii) is the \emph{parity axiom}. (Dropping parity amounts to dropping (ii) and
relaxing (\ref{eq:star=}) from  $C_1(\phi,\psi)=\thalf\{\phi,\psi\}$  to
$C_1(\phi,\psi)-C_1(\psi,\phi)=\{\phi,\psi\}$.)
In axiom (iii), we have extended pointwise complex conjugation  to
$\A[t]$ so that  $\ovl[\phi t^i]=\ophi t^i$.
Axiom (iv) means that if $\phi$ or $\psi$ vanishes identically on $T^*U$,
where $U$ is open in $\XR$,
 then $C_p(\phi,\psi)$ vanishes identically on $T^*U$.
Then the operators $C_p=C_p(\cdot,\cdot)$, and hence the star product,
extend  naturally  from $\PolTXR$ to $\Pol(T^*U)$. 

Axiom (v) is often called \emph{strong invariance} -- we use the term ``equivariant".
This is an   important notion because it corresponds to equivariant quantization of
symbols (see \S\ref{sec_bqsm}). Strong invariance implies the weaker notion of 
\emph{invariance}, which means that the operators $C_p$ are $\GR$-invariant. 

We note that   (iv)  is   much weaker  than
the familiar axiom that requires  locality on  $\TXR$. Indeed  
locality means that the   $C_p$, and hence the star product, are
bidifferential.  It turns out that bidifferentiality is too strong a geometric requirement
in our situation, but   we believe  it can be modified in a controlled way consistent 
with  (iv);  see \S\ref{sec_Lax}.

At $t=1$, $\star$  specializes   to a noncommutative product on
$\B=\A[t]/(t-1)$; this works  because of axiom  (i).
Then $\B$ has an  increasing algebra filtration (defined by 
the grading on $\A$) and the obvious vector space isomorphism $\Q:\A\to\B$
induces a graded Poisson algebra  isomorphism  from $\A$ to
$\gr\B$. Via $\Q$, the structures on $\A$ pass over to $\B$. Axiom 
(ii)  implies that $\al$ defines an algebra anti-involution $\be$ on $\B$
and (iii) implies that
complex conjugation on $\A$ defines an anti-linear  algebra
involution  $a\mapsto\oa$ on $\B$.
By (v),   $\B$ acquires a \rep\ of $\GR$ compatible with everything.

We can find a nice  candidate for $\B$  by asking for  compatibility between 
deformation quantization and geometric quantization.
Geometric quantization  of $T^*\XR$ produces  the Hilbert space
$\cH$ of square integrable  half-densities on $\XR$, where
\emph{half-densities} are the (complex-valued) sections of the  half-density line
bundle $\Ehalf_{\XR}$. If we ask that $\B$ operates on (a dense subspace of) 
$\cH$,   the obvious candidate for $\B$ is the  algebra $\DXR$ of smooth 
differential  operators for   $\Ehalf_{\XR}$.

Fortunately,  $\DXR$ already has all the structure discussed above.  
It has the order filtration and the principal symbol map identifies $\gr\DXR$ with
$\A$. Also $\DXR$ admits the  pointwise complex conjugation  map
$a\mapsto\oa$  defined by  $\oa(\sig)=\ovl[{a(\ovl[\sig])}]$ where 
$\sig\mapsto\ovl[\sig]$ is pointwise complex conjugation  of half-densities.
There is a canonical  $\GR$-invariant algebra anti-involution $\be$
of $\DXR$ such that $\be(\phi)=\phi$ for $\phi\in\A^0$ and
$\be(\etah)=-\etah$ if  $\etah$ is the Lie derivative of a  vector  field $\eta$ on
$\XR$.  Then $\be$ induces $\al$ upon taking principal symbols.
Finally, we have a compatible   \rep\ of $\GR$  on $\DXR$
because  the line bundle $\Ehalf_{\XR}$ is $\GR$-homogeneous.

We have the  Lie algebra homomorphism
$\g\to\DXR$, $x\mapsto\etah^x$, where $\eta^x$ is the complex vector field  on
$\XR$ defined by $x$.
The  \rep\ of $\g$ on $\DXR$  by the operators
$[\etah^x,\cdot]$ corresponds to the natural  \rep\ of $\GR$.
(As is often done, we   complexify the group \rep\ at the Lie algebra level.)

\section{Quantizing symbols into differential operators equivariantly}  
\label{sec_bqsm}  
Now that we have decided upon $\B=\DXR$, 
we can reformulate  our  star product problem in terms of quantization maps.
To begin with,  we can axiomatize the properties of 
our vector space isomorphism $\Q:\A\to\B$ from \S\ref{sec_eqprob}:  
\begin{equation}\label{eq:bqsm_list} 
  \begin{array}{ll}
     \rm{(i)}&\mbox{if $\phi\in\A^d$ then the principal symbol of
                       $\Q(\phi)$ is $\phi$}\\[1pt]
     \rm{(ii)}   &\mbox{$\Q(\phi\sal)=\Q(\phi)\sbe$}  \\[1pt]
     \rm{(iii)} &\mbox{$\Q(\ophi)=\ovl[\Q(\phi)]$}  \\[1pt]
     \rm{(iv)} & \Q  \mbox{ is local on } \XR   \\[1pt]
     \rm{(v)} & \Q(\mu^x)=\etah^x \mbox{ and }
                       \Q(\{\mu^x,\phi\})= [\etah^x,\Q(\phi)]
                   \,\,\mbox{ if } x\in\g\\[2pt]
  \end{array} 
\end{equation}
Axiom (iv) means that if $U$ is open in $\XR$ and $\phi$   
vanishes identically on $T^*U$ then the differential operator  $\Q(\phi)$
vanishes identically on  $U$. In (v), we used the semisimplicity of $\g$ to get
$\Q(\mu^x)=\etah^x$.
Axiom (v)  means that $\Q$ is $\g$-equivariant. This amounts to
$\GR$-equivariance.

We call $\Q$ a \emph{$\GR$-equivariant quantization map}.
We can recover $\star$ from $\Q$  by the formula
$\phi\star\psi=\Q_t\inv(\Q_t(\phi)\Q_t(\psi))$ where   
$\Q_t(\phi t^p)=\Q_t(\phi)t^{j+p}$ if $\phi\in\A^j$.
In this way, we get a bijection between graded equivariant star products on $\A$ 
and  equivariant quantization maps (up to algebra automorphisms of $\DXR$  which
are compatible with principal symbols, the $\GR$-action, etc.).

\section{Algebraicity of the EDQ problem for $\TXR$}  
\label{sec_RR}  

Since $X$ is a complex algebraic (projective) variety, we can consider  the
\emph{algebraic parts} $\RR$ and $\DD$ of $\PolTXR$ and $\DXR$. 
By this we mean that $\RR$ is the subalgebra of $\PolTXR$ 
corresponding,  by restriction of functions, to the algebra
$\RTX$ of regular functions   on the quasi-projective  variety $T^*X$.
Similarly, $\DD$ is the subalgebra of $\DXR$ corresponding to the algebra 
$\DX$ of  algebraic twisted differential operators  for the
(locally defined)  square root of the canonical bundle.
We have  $\RTX\simeq\RR$ and $\DX\simeq\DD$; this follows since  
$\TXR$ is a real form of $\TX$.
 
The action of $G$ on $X$ induces natural \reps\ of $G$ on $\RTX$ and $\DX$, and
hence on  $\RR$ and $\DD$, which are both locally finite and completely reducible.
(Locally finite for $G$ means that every vector lies in a finite-dimensional $G$-stable
subspace.) 
Thus $\RR$ and $\DD$ have    more symmetry than $\A$ and $\B$. 

Clearly $\RR$ contains the algebra  $\RRm$ generated
by the momentum   functions $\mu^x$, $x\in\g$, and $\DD$ contains the algebra
$\DD_{\eta}$ generated by the twisted vector fields  $\etah^x$.
Soon (\S\ref{sec_thm} onwards)
we will restrict to the case where $\RR=\RRm$ and $\DD=\DD_{\eta}$.

We can formulate the notion of a  \emph{graded $\GR$-equivariant star product on 
$\RR$}    using the same axioms  as in \S\ref{sec_eqprob}.
(Axiom  (\ref{eq:star_list})(iii)  makes sense because $\RR$ is stable under
$\phi\mapsto\ophi$. Axiom (\ref{eq:star_list})(iv) is  vacuous as stated.)
Similarly, once we establish $\gr\DD=\RR$,
we can  formulate the notion of 
a \emph{$\GR$-equivariant quantization map $\bq:\RR\to\DD$} using the same
axioms   as  in   \S\ref{sec_bqsm}.  
In both cases, $\GR$-equivariance easily implies $G$-equivariance.

\begin{prop}\label{prop:res} 
\begin{itemize}
\item[(i)]  Any graded $\GR$-equivariant  star product  
$\star$  on  $\A$ restricts to a graded $G$-equivariant star product  on $\RR$.
\item[(ii)]  We have $\gr\DD=\RR$.
Any $\GR$-equivariant quantization map $\Q:\A\to\B$ restricts to a
$G$-equivariant quantization map $\bq:\RR\to\DD$.
\end{itemize}
\end{prop}   

\begin{proof} 
(i) We just need to show that  $\star$ restricts to $\RR$, i.e.,  if $\phi$ and 
$\psi$ belong to  $\RR$ then $\phi\star\psi$ belongs  to  $\RR[t]$.
 $\GR$-invariance of   $\star$   implies that $\star$ restricts to
the $\GR$-finite part $\Afin$ of $\A$. ($\Afin$ is the subalgebra consisting of 
functions which  lie  in a finite-dimensional $\GR$-stable subspace of $\A$.)
We will show that $\RR=\Afin$. 
Certainly  $\RR$ lies in $\Afin$ since $G$ is locally finite on $\RR$. 

Now $\RTX[d]$ identifies with the space $\cP^d=\Polhol[d]$
of holomorphic functions on $T^*X$ which are homogeneous 
degree $d$ polynomials on the cotangent fibers.
This follows from compactness of $X$. We will show that
any $\phi$  in $\Afin^d$ extends 
to a holomorphic function $\phi\shol$ in $\cP^d$.

Now   $\cP^d$   is the (finite-dimensional) space of holomorphic sections of a 
finite rank vector  bundle $E^d$ over  $X$ and    
$\A^d$ is the space of smooth sections of $E^d|_{\XR}$.   
Both  $E^d|_{\XR}=\GR\times_{\PR} V^d$ and $E^d=G\times_{P} V^d$ 
are homogeneous vector bundles where $V^d$ is the base fiber  of $E^d$. 
Consequently we can identify 
$\A^d$ and  $\cP^d$ with  certain spaces    
of  functions  $\GR\to V^d$ and $G\to V^d$ in the familiar way.
If $G$ is simply connected,
every smooth $\GR$-finite function  on $\GR$   extends uniquely 
to a holomorphic  ($G$-finite) function   on $G$. Using this, it follows easily that
$\phi$ extends to $\phi\shol$ as desired.

(ii) We have  $\gr\DX=\RTX$ by \cite[Lem. 1.4]{Bo-Br} --  
their  result goes through to the twisted case with the same proof.
So $\gr\DD=\RR$.  

Now $\DD$ lies in the $\GR$-finite part $\Bfin$
of  $\B$,  since $G$ is locally finite on $\DD$.  
Using principal symbols, we find $\DD=\Bfin$,
since $\gr\Bfin$ lies in
$\Afin=\RR$. (In fact this proves $\gr\Bfin=\Afin$.)
Clearly $\Q(\Afin)=\Bfin$, and so $\Q(\RR)=\DD$. 
\end{proof}

We regard $\RR$ and $\DD$ as algebraic   models  of  $\A$ and $\B$.  
We  know that  $\RR$ and $\DD$ are finitely generated algebras, and $\RR$ is
finite as a module over $\RRm$ (see e.g. \cite{B-K:jams}). 
At first sight $\RR$  may seem to be too small to encode enough information about
$\A$. For instance, $\RR^0=\C$  while $\A^0$ is the  infinite-dimensional
algebra of smooth complex valued functions on $\XR$.  
But already $\RRm$ is ``big  enough" in the sense that 
\begin{lem}\label{lem:diff_on_A} 
Any differential operator  $A$ on $\TXR$ is uniquely determined by
its values  $A(\phi)$ where $\phi$ belongs to $\RRm$.
\end{lem} 
\begin{proof} 
This follows since the momentum functions $\mu^x$, $x\in\gR$, form a complete set 
of functions  (i.e. their differentials span the cotangent spaces) over some open dense
set $W$ in $\TXR$. Indeed, the image of  the   moment map
$\muR:\TXR\to\gR^*$  is the closure of a single
nilpotent orbit $\OR$, and we can choose  $W=\muR\inv(\OR)$.
\end{proof}

Proposition \ref{prop:res}
suggests that we might try to solve to the EDQ problem posed in 
\S\ref{sec_eqprob} by finding a ``preferred", or particularly natural,
solution to the analogous   $G$-equivariant problem for $\RR$ and $\DD$,
and then trying to extend that solution to   $\A$ and $\B$. 
This extension problem would have a unique solution, on account of 
Lemma \ref{lem:diff_on_A}, if we found a star product on $\RR$ given  by
bidifferential operators on $\TX$. This same kind of argument can still work
if  bidifferentiality is violated in a controlled way, by inverting certain
nice invertible  operators; 
see \cite[Th. 5.1]{me:RPn} for an example and  \S\ref{sec_Lax} for a conjecture.

In the next section we find a   preferred  $G$-equivariant graded
star product on $\RR$. We do this under the hypothesis that $\RR=\RRm$.
This is a  hypothesis on $(G,X)$ which is
satisfied for instance if (i) $G=SL_n(\C)$ and  $X$ is arbitrary (\cite{K-P}), 
or   (ii) $G$ is arbitrary but $X$ is the full flag manifold.

This hypothesis was important in \cite{Bo-Br} in studying noncommutative analogs
of $\RTX$; it is  equivalent  (\cite[Th. 5.6]{Bo-Br}) to the condition that
the holomorphic moment map $\mu:\TX\to\g^*$ has  \emph{good geometry}
in the sense that $\mu$  is generically $1$-to-$1$  and its image 
in $\g^*$ is a normal variety.  These conditions have been studied a lot in geometric
\rep\  theory, especially since the image of $\mu$ is  the closure of a single  
nilpotent coadjoint orbit $\OO$ of $G$.

Suppose $\RR=\RRm$. Then  
$\RR =\Sg/I$ where $I$ is the (graded) ideal of  functions in $\Sg$ 
which  vanish  on $\OO$, and $\DD=\Ug/J$ where $J$ is a $2$-sided ideal
in  $\Ug$ with $\gr J=I$. The ideal $I$ contains all casimirs (i.e., $G$-invariants in
$\oplus_{d=1}^\infty\Sg[d]$). The casimirs generate $I$  if and  only if $X$ is the
full flag variety.

\section{A preferred star product on $\RR$}  
\label{sec_thm}
Suppose $\phi\star\psi$ is a  graded $G$-equivariant
star product  on  $\RR$ (see \S\ref{sec_RR}).  This defines  a noncommutative
associative product $\circ$ on $\RR$ where
$\phi\circ\psi$ is the specialization at $t=1$ of $\phi\star\psi$.  Then we obtain a
\rep\ $\pi$ of $\gog$ on $\RR$ given by  
$\pi^{x,y}(\phi)=\mu^x\circ\phi-\phi\circ\mu^y$.
Notice that the equivariance axiom \textup{(\ref{eq:star_list})(v)} says
that the  quantum  operator $\pi^{x,x}$ coincides with the  classical operator
$\{\mu^x,\cdot\}$.

\begin{thm}\label{thm:main} 
Assume $\RR$ is generated by  $\mu^x$, $x\in\g$.
Suppose $\star$ is a  graded $G$-equivariant star product  on $\RR$ where 
$\star$ corresponds to a $G$-equivariant quantization map
$\bq:\RR\to\cD$. \textup{(}Such maps $\bq$ always exist.\textup{)} Then
\begin{itemize}   
\item[(I)] The \rep\  $\pi$ of $\gog$ on $\RR$  is irreducible and unitarizable, i.e.,
there exists a unique positive definite invariant hermitian form $\apair{\cdot}{\cdot}$
on
$\RR$ with  $\apair{1}{1}=1$.
\item[(II)] There is a unique choice of $\bq$, and hence a unique choice of $\star$, 
such that the grading $\RR=\oplus_{d=0}^\infty\RR^d$ is orthogonal with respect 
to  $\apair{\cdot}{\cdot}$. Then  
\begin{equation}\label{eq:pixyphi=} 
\pi^{x,y}(\phi)=\mu^{x-y}\phi+\thalf\{\mu^{x+y},\phi\}+\La^{x-y}(\phi)
\end{equation}
where $\La^x$, $x\in\g$, are  certain operators on $\RR$. 
\end{itemize} 
\end{thm}
\begin{proof} 
The proof occupies \S\ref{sec_exist}--\ref{sec_II}.
\end{proof}

We now discuss what unitarizable means and introduce some notations.
To begin with, the restriction of $\pi$ to $\g^{diag}=\{(x,x)\,|\, x\in\g\}$, i.e,  
the $\g$-\rep\  on $\RR$  given  by the operators $\pi^{x,x}$,  
corresponds to the natural  $G$-\rep\  on $\RR$. 
Thus $\RR$ is a $(\gog,G)$-module in the sense of   Harish-Chandra.

Now  \emph{unitarizability} of  $\pi$ means that there is a positive definite hermitian
inner product $\apair{\cdot}{\cdot}$ on $\RR$ which is invariant for
$\gs=\{(x,\sig(x))\,|\, x\in\g\}$, i.e., the operators 
$\pi^{x,\sig(x)}$ are skew-hermitian. Here $\sig$ is a  fixed Cartan involution of
$\g$; we   choose $\sig$ compatible with $\GR$ so that $\sig$ extends
a Cartan involution of  $\gR$.  Then $\sig$ corresponds to a maximal compact 
subgroup $\Gc$ with Lie algebra $\gc=\{x\in\g\,|\, x=\sig(x)\}$.   
E.g., if $\gR=\fsl_n(\R)$, then  take  $\sig(x)=-\ox^t$ so that $\gc=\fsu_n$.

By a theorem of Harish-Chandra, the  operators $\pi^{x,\sig(x)}$ then correspond  to a
unitary representation of  $G$ on the  Hilbert space completion  $\hatRR$ of $\RR$
with respect to $\apair{\cdot}{\cdot}$. If the   $\RR^d$ are orthogonal, then  
 $\hatRR$ is the Hilbert space direct sum 
$\widehat{\oplus}_{d=0}^\infty\RR^d$.
Notice that we end up with two very different actions of $G$: the graded algebraic 
action on $\RR$ corresponding to $\g^{diag}$ and the 
unitary action on $\hatRR$ corresponding to $\gs$.

\section{Existence proof for    $\bq$}  
\label{sec_exist}   

A $G$-equivariant  quantization map $\bq$ is completely determined by
the subspaces  $\F^d=\bq(\RR^d)$.
This is immediate from (\ref{eq:bqsm_list})(i).
Then the   decomposition  $\DD=\oplus_{d=0}^\infty\,\, \F^d$
``splits the order filtration"  in the sense that  
$\oplus_{d=0}^p\,\F^d=\DD_{\le p}$.
Referring to (\ref{eq:bqsm_list}) again, we see that
the spaces $\F^d$ are stable under $\be$, complex conjugation, and $\g$ (which 
acts  by  $A\mapsto [\etah^x,A]$). Conversely, any such splitting  
corresponds to a choice of $\bq$.  

\begin{lem}\label{lem:bqexists}   
We can always construct a $G$-equivariant  quantization map $\bq:\RR\to\DD$. 
If the \rep\ of $G$ on $\RR$ is multiplicity free,   there is  
only one  choice for $\bq$.
\end{lem}  
\begin{proof} 
By complete  reducibility, we can  find  a $\g$-stable complement $\cG^d$   to 
$\DD_{\le d-1}$ inside $\DD_{\le d}$.
This gives a $\g$-stable splitting of the order filtration; let $\bp$ be the 
corresponding  quantization map. The spaces $\cG^d$ may fail to be stable
under $\be$ and/or complex conjugation.  To remedy this, we  ``correct" $\bp$ by
putting $\bp'(\phi)=\half\left(\bp(\phi)+\bp(\phi\sal)\sbe\right)$ and then
$\bp''(\phi)=\half\left(\bp'(\phi)+\ovl[\bp'(\ophi)]\,\right) $.
Now $\bp''$  is  a valid choice for $\bq$.

If $\RR$ is multiplicity free, then  $\cG^d$ is unique
for each $d$, and so $\bp$ is the unique choice for $\bq$.
Notice that uniqueness of $\bq$ does not require (ii)-(iv) in
(\ref{eq:bqsm_list}).
\end{proof} 

In the multiplicity free case,  the method explained in Remark \ref{rem:symm}
gives a sort of formula for  $\bq$.  We note that 
$\RR$ is multiplicity free whenever the parabolic subgroup $P$
(where $X=G/P$) has the property that its unipotent radical is abelian. 
For $G=SL_n(\C)$,   this happens when  $X$ is a grassmannian. The
full classification of multiplicity free cases  is well known.

In general, there will be infinitely many choices for $\bq$; we can show using
filtration splittings that  the set of   choices has the structure of an  infinite 
dimensional affine space. 

\section{Proof of (I) in Theorem \ref{thm:main}}  
\label{sec_I}   
The quantization map $\bq$ intertwines our \rep\ $\pi$ of $\gog$ on $\RR$ with
the  \rep\ $\Pi$ of $\gog$ on $\DD$ given by  
$\Pi^{x,y}(A)=\etah^xA-A\etah^y$. Indeed, 
$\bq(\phi\circ\psi)=\bq(\phi)\bq(\psi)$ and so
$\bq(\pi^{x,y}(\phi))=\Pi^{x,y}(\bq(\phi))$. 

Therefore  proving $\pi$ is irreducible and unitarizable  reduces to  proving  
$\Pi$ is irreducible and unitarizable. For this, we need our hypothesis
that $\RR$ is generated by the $\mu^x$; we use this freely from now on.

Each  $\eta^x$  extends uniquely to a holomorphic vector field 
$\xi^x$ on  $X$.    We can regard $X$ as a real manifold. Then 
$\DXsmooth$ is the algebra of smooth differential
operators on the space $\GEhalf$  of smooth half-densities on $X$.
The  Lie derivative $\xih^x$ lies in $\DXsmooth$.
  
The map $\etah^{x}\mapsto\xih^{x}$ extends naturally to   a
$G$-equivariant algebra embedding $A\mapsto\Ahol$ of
$\DD$ into $\DXsmooth$.
Let $\delta$ be the unique  $\Gc$-invariant positive real density  on $X$ such that
$\int_X\delta=1$. Let $\dhalf$ be the  positive square root of $\delta$.
We map $\DD$ into $\GEhalf$ by
\begin{equation}\label{eq:Delta} 
\Delta(A)=\Ahol(\dhalf)
\end{equation}
Now   $\DD$ acquires  the $\Gc$-invariant  hermitian pairing    
$\ga(A,B)=\int_X \Delta(A)\,\ovl[\Delta(B)]$.

\begin{prop}\label{prop:ga} 
$\ga$ is  $\gs$-invariant and positive definite. 
\end{prop}
\begin{proof}
$\gs$-invariance means that the operators $\Pi^{x,\sig(x)}$ are skew-hermitian, or
equivalently, the adjoint of  $\Pi^{x,0}$ is $-\Pi^{0,\sig(x)}$. So we want to show
\begin{equation}\label{eq:ga=ga} 
\ga(\etah^xA,B)=\ga(A,B\etah^{\sig(x)})
\end{equation}
We have $\ga(\etah^xA,B)=\int_X(\xih^x\Delta(A))\,\ovl[\Delta(B)]
=-\int_X\Delta(A)\,(\xih^x\ovl[\Delta(B)])$; 
the last equality holds because $\int_X\xih^x(\al\be)=0$ for any
half-densities $\al,\be$.  

$\Gc$-invariance of $\dhalf$ means  that 
$\xih^x+\ovl[\xih^x]$  kills $\dhalf$ if   $x\in\gc$,  or equivalently 
$\xih^x+\ovl[\xih^{\sig(x)}]$  kills $\dhalf$ if   $x\in\g$.   
Using this and  the commutativity of  holomorphic and anti-holomorphic operators   
we find $\xih^x\,\ovl[\Delta(B)]=
-\ovl[B\shol \xih^{\sig(x)}](\dhalf)
=-\ovl[\Delta(B\etah^{\sig(x)})]$
and so we get  (\ref{eq:ga=ga}).

For positive definiteness, we just need to show that   $\Delta$ is $1$-to-$1$ on
$\DD$.  We expect there is  a geometric proof of this, but we have not worked that
out. Instead, we will use   results from representation theory. This argument will be
clear for experts in these matters  and too technical for everyone else; so we   just
sketch it briefly.

$\Delta$ is $\Gc$-equivariant and so $\Delta$ maps 
$\DD$ into $\GEhalffin$. 
We have $\DD=\Ug/J$ where    $\gr J=I$.
We can show that $\Delta$  is the same   (up to scaling)
as the map  from $\Ug/J$ to $\GEhalffin$ defined by Conze-Berline and Duflo in 
\cite[\S5.3 and Cor. 6.3]{CB-D}. (This is the ``$\pi=0$" case in their notation.)
This follows  easily because both maps are 
$\gog$-equivariant; here $\gog$ acts on  $\GEhalf$ by 
the twisted vector fields  $\xih^{x,y}=\xih^x+\ovl[\xih^{\sig(y)}]$. 
Next we need a suitable criterion for injectivity of the 
Conze-Berline--Duflo map. We find
it in Vogan's result \cite[Prop. 8.5]{Vo1984} on  injectivity of certain maps of 
induced modules into produced modules;  see also \cite[\S6]{Vo1990}. 
\end{proof}

\begin{cor}\label{cor:HC} 
$\DD$ is isomorphic, via  $\Delta$, to the Harish-Chandra module   
of the natural unitary \rep\  of $G$ on $L^2(X,\Ehalf_X)$.
\end{cor} 
\begin{proof} 
The Harish-Chandra module is $\GEhalffin$.
We just established injectivity of $\Delta$.  Surjectivity  
follows by \cite[Prop. 5.5]{CB-D} or by  directly checking that
the source and target contain the same irreducible $\Gc$-representations with the
same multiplicities. 
\end{proof}

\begin{cor}\label{cor:ga} 
We have $\ga(A,B)=\ga(B\ssig A,1)$ where $B\mapsto B\ssig$ is an anti-linear
algebra involution of $\DD$. 
\end{cor}
\begin{proof} 
Since the $\etah^x$ generate $\DD$, we may   assume
$B=\etah^{y_1}\cdots\etah^{y_m}$. Then (\ref{eq:ga=ga}) gives 
$\ga(A,B)=\ga(B\ssig A,1)$ where
$B\ssig=\etah^{\sig(y_1)}\cdots\etah^{\sig(y_m)}$.
This map $B\mapsto B\ssig$ gives a well-defined anti-linear algebra involution of
$\DD$; indeed we have $\ovl[\Delta(B\sbe)]=B\ssig(\dhalf)$.
\end{proof}

The formula $\T(A)=\ga(A,1)$ defines a linear functional $\T$ on $\DD$. Explicitly,
\begin{equation}\label{eq:T} 
\T(A)=\int_X A\shol(\dhalf)\dhalf
\end{equation}
\begin{prop}\label{prop:T} 
$\T$ is the unique $\Gc$-invariant  linear  functional on $\DD$ with $\T(1)=1$.
Moreover $\T$ is a trace.  We have 
$\ga(A,B)=\T(AB\ssig)$.
\end{prop}
 
\begin{proof} 
Clearly $\T$ is $\Gc$-invariant. Then
$\T:\DD\to\C$ is the unique invariant linear projection  because 
the $\Gc$-action on $\DD$ is
completely reducible and  the constants  are the only $\Gc$-invariants in $\DD$
(since  the constants are the only $\Gc$-invariants in $\RR$).

$\T$ is   $\g$-invariant, i.e.,
$\T([\etah^y,A])=0$. We  write this as $\T(\etah^yA)=\T(A\etah^y)$. Iteration 
gives
$\T(\etah^{y_1}\cdots\etah^{y_k}A)=\T(A\etah^{y_1}\cdots\etah^{y_k})$.
This proves $\T(BA)=\T(AB)$ . 
\end{proof}

Now we can show that
$\ga$ is the unique $\gs$-invariant hermitian form on $\DD$ such that
$\ga(1,1)=1$. Indeed suppose $\nu$ is any such  form. Then
$\nu(A,1)=\T(A)$ by the uniqueness of $\T$. So (\ref{eq:ga=ga}) gives 
$\nu(A,B)=\nu(B\ssig A,1)=\T(B\ssig A)=\ga(A,B)$. 
This uniqueness of $\ga$ implies $\Pi$ is irreducible.

This completes the proof of    Theorem \ref{thm:main}(I).  
Once $\bq$ is chosen,  $\apaircd$ is given by
\begin{equation}\label{eq:apair=} 
\apair{\phi}{\psi}=\ga(\bq(\phi),\bq(\psi))=
\T\left(\bq(\phi)\bq(\psi)\ssig\right)
\end{equation}

Finally we note that  the irreducibility of $\Pi$ implies (and vice versa) 
\begin{cor}\label{cor:simple} 
$\DD$ is a simple ring.
\end{cor}

\section{Proof of (II) in Theorem \ref{thm:main}}
\label{sec_II}   

The graded pieces $\RR^d$ 
are orthogonal with respect to $\apair{\cdot}{\cdot}$ iff
their images $\bq(\RR^d)$ are orthogonal with respect to $\ga$.
So we have only one possible choice of $\bq$, namely the one such that 
$\bq(\RR^d)=\cV^d$ where
$\oplus_{d=0}^\infty\cV^d$ is the   $\ga$-orthogonal splitting of the  order
filtration of $\DD$. 
According to \S\ref{sec_exist},  we need to check
\begin{lem}\label{lem:cV} 
$\cV^d$ is stable under $\be$, complex conjugation, and $\g$.
\end{lem}
\begin{proof} 
$\Gc$-invariance of $\ga$ implies that $\cV^d$ is $\g$-stable.
The other two follow easily   using  $\ga(A,B)=\T(AB\ssig)$ and   the properties
$\T(A\sbe)=\T(A)$, $\T(\oA)=\ovl[\T(A)]$.
\end{proof}
Thus $\oplus_{d=0}^\infty\cV^d$ defines $\bq$. Then $\bq$ defines a graded
$G$-equivariant star product $\star$ on $\RR$; this is the only one for  which the
sum $\oplus_{d=0}^\infty\RR^d$ is $\apaircd$-orthogonal.
\begin{prop}\label{prop:RjRk} 
This  star product $\star$ satisfies 
\begin{equation}\label{eq:RjRk} 
\RR^j\star\RR^k
                    \subseteq\RR^{j+k}\oplus\cdots\oplus\RR^{|j-k|}t^{2\min{(j,k)}}
\end{equation}
\end{prop}
\begin{proof} 
Since $\star$ is graded, it suffices to consider $\circ$.
Let $\ell(\phi)$ and $r(\phi)$ denote respectively left and right 
$\circ$-multiplication by  $\phi$.  
It is easy to check that the map $\mu^x\mapsto\mu^{\sig(x)}$ extends to a
graded anti-linear algebra   involution 
$\phi\mapsto\phi^{\sig(x)}$ of $\RR$; this follows
because the nilpotent orbit $\OO$ is $\sig$-stable.
It follows by
(\ref{eq:ga=ga})   that the adjoint of  $\ell(\phi)$ is $r(\phi\ssig)$.

Suppose $\phi\in\RR^j$. Then the  highest degree term in $\phi\circ\psi$,
namely   $\phi\psi$,  occurs in degree  $j+k$. If $\nu\in\RR^d$ and $\nu$ occurs in 
$\phi\circ\psi$, then $\psi$ occurs in $\nu\circ\phi\ssig$ and so
$d+j\ge k$. Similarly  $d+k\ge j$. So $d\ge |j-k|$.
\end{proof}
\begin{cor}\label{cor:3term}   
For $x\in\g$ and $\phi\in\RR$ we have
\begin{equation}\label{eq:mux*phi} 
\mu^x\star\phi=\mu^x\phi+\half\{\mu^x,\phi\}t+\La^x(\phi)t^2
\end{equation}
where $\La^x$ is the adjoint with respect to $\apair{\cdot}{\cdot}$ of ordinary
multiplication by $\mu^{\sig(x)}$.
\end{cor}
\begin{proof} 
Certainly  (\ref{eq:RjRk}) implies   (\ref{eq:mux*phi}) where
$\La^x(\psi)=C_2(\mu^x,\psi)=C_2(\psi,\mu^x)$.
Now suppose $\phi\in\RR^j$ and $\psi\in\RR^{j+1}$. Because of orthogonality of the
spaces $\RR^d$ we find  
$\apair{\phi}{\La^x(\psi)}=\apair{\phi}{\psi\circ\mu^x}=
\apair{\,\mu^{\sig(x)}\circ\phi}{\psi}=\apair{\,\mu^{\sig(x)}\phi}{\psi}$.
\end{proof}
Now (\ref{eq:mux*phi}) gives (\ref{eq:pixyphi=}). This concludes the proof of
Theorem 
\ref{thm:main}.

\begin{remark}\label{rem:symm}  
We know  another method for constructing a  $G$-equivariant quantization map 
$\bfr:\RR\to\DD$.   We start with the positive
definite hermitian  pairing $\bh(f,g)=\del_g(f)$ on $\Sg$, where $\del_x$ is the
constant coefficient vector field  on $\g$ defined by
$\del_x(y)=-(\sig(x),y)$ and $\del_{g_1g_2}=\del_{g_1}\del_{g_2}$.
Let   $H$   be the $\bh$-orthogonal complement to $I$ in $\Sg$ where
$\RR=\Sg/I$.  Then $H=\oplus_{d=0}^\infty H^d$ is  graded.   
We put $\F^d=p(s(H^d))$,  
where $s:\Sg\to\Ug$  is the symmetrization map and $p:\Ug\to\DD$ is the natural
projection.  
Then $\DD=\oplus_{d=0}^\infty\F^d$ is a $\g$-stable splitting of the order filtration.
It is easy to see that $\F^d$ is  stable under $\be$, complex conjugation, 
and also $\sig$. So  by \S\ref{sec_exist}   this splitting defines $\bfr$.

Here is a formula for  $\bfr$:  if we pick a basis $x_1,\dots,x_m$ of $\g$ and
$\sum a_{i_1,\dots,i_d}\,x_{i_1}\cdots x_{i_d}$ lies in $H^d$,  then  
\begin{equation}\label{eq:bq_symm} 
\bfr\left(\tsum a_{i_1,\dots,i_d}\,\mu^{x_{i_1}}\cdots\mu^{x_{i_d}}\right)
=\frac{1}{d!}\tsum_{\tau} 
\,a_{i_1,\dots,i_d}\,\etah^{x_{i_{\tau(1)}}}\cdots\etah^{x_{i_{\tau(d)}}}
\end{equation}
where we sum over all permutations $\tau$ of $\{1,\dots,d\}$.
 
We conjecture that $\F^d=\cV^d$, or equivalently, that $\bfr=\bq$.
This is obviously true in the multiplicity free case by
uniqueness (Lemma \ref{lem:bqexists}).  Analytic methods
may well be needed to show  $\bfr=\bq$, just as  we needed
integration to establish the  positivity of $\ga$ (or even the weaker fact that
$\ga$ is non-degenerate on each space $\D_{\le d}$).
 
Suppose  $X$ is the full flag manifold. Then  $H$ is Kostant's  space of harmonic
polynomials, and  $\bfr$    is  simply a 
$\rho$-shifted version of the map   constructed by Cahen and Gutt in \cite{C-G} for   
the principal nilpotent orbit case.  
\end{remark} 

\section{The operators $\La^x$ on $\RR$}  
\label{sec_Lax}   

In Corollary \textup{\ref{cor:3term}} we saw that our star product $\star$
produces   operators $\La^x$, $x\in\g$, on $\RR$.
Conversely, the  $\La^x$ completely determine $\star$.
This follows because,  if we know the $\La^x$, then  using associativity we can
compute  $(\mu^{x_1}\cdots\mu^{x_n})\star\psi$   by
induction on $n$. Here (\ref{eq:mux*phi}) provides the first step
$n=1$, and also it propels the   induction.

Several nice properties follow from Corollary \textup{\ref{cor:3term}}:  
\begin{itemize}
\item[\rm(i)] $\La^x$ is graded of degree $-1$, i.e., i.e.,
    $\La^x(\RR^j)\subseteq\RR^{j-1}$.
\item[\rm(ii)] The $\La^x$ commute and generate a graded
    subalgebra of $\End\,\RR$  isomorphic to $\RR$.
\item[\rm(iii)]  The $\La^x$ transform in the
    adjoint \rep\ of $\g$, i.e., $[\Phi^x,\La^y]=\La^{[x,y]}$
    where $\Phi^x=\{\mu^x,\cdot\}$
\item[\rm(iv)]  The map $\g\times\g\to\C$, $(x,y)\mapsto\La^x(\mu^y)$,  is a
         non-degenerate $\g$-invariant symmetric complex bilinear pairing.
\end{itemize}

The $\La^x$ are  not  differential operators on $\RR$ in general.
Indeed differentiality   fails when $\GR=SL_{n+1}(\R)$ and $\XR=\RP^n$. 
In that case  $\La^x$ is a reasonably
nice operator as it is the left quotient of an algebraic  
differential operator  $L^x$ (of order $4$) on the closure on $\OO$  by   the 
invertible  operator $(E+\half[n]) (E+\half[n]+1)$. Moreover $L^x$  extends 
to $T^*\CP^n$.   See \cite{A-B:starmin} and  \cite{me:RPn}.

The $\La^x$ determine $\star$ in a rather simple way, and so their
failure  to be differential should control the failure of $\star$ to be
bidifferential. This issue is important in understanding if and how $\star$ extends 
from  $\RR$ to  $\A$ (see the discussion in \S\ref{sec_RR}). 

We conjecture that $\La^x$ is of the form $\La^x=P\inv L^x$ where 
(i) $P$ and $L^x$ are  algebraic  differential operators on $\TX$,
(ii)  $P$ is  $G$-invariant and   \emph{vertical} so that  $P$ ``acts along the fibers 
of  $\TX\to X$"  (iii)  $P$ is invertible on $\RR$, in fact $P$ is diagonalizable with  
positive spectrum and (iv)   the formal order of  $P\inv L^x$   is $2$. 
 
Motivated by this conjecture, J-L. Brylinski and the author construct in
\cite{us:vertical} all invariant
vertical differential operators on cotangent bundles of grassmannians.  

The $\XR=\RP^n$ case discussed above is an example where (i)-(iv) works.
That example was part of a quantization program for  minimal nilpotent orbits
(see \cite[\S1]{A-B:exotic}). In fact, our conjecture here  arises   from
a larger program  we have on quantization of general nilpotent orbits.
A proof our conjecture,  coming most likely out of  properties of  $\apaircd$,
would  give more evidence for our program.

\section{The inner product $\apaircd$ on $\RR$}  
\label{sec_inner}   

In Theorem \ref{thm:main},   the hermitian form $\apaircd$ completely 
determines  the  star product $\star$, and vice versa. To show this, it suffices 
(see \S\ref{sec_Lax}) to show that knowing  $\apaircd$ is equivalent to knowing 
the $\La^x$. Certainly $\apaircd$ produces  $\La^x$, as $\La^x$ is (Corollary
\ref{cor:3term})  the adjoint of  $\phi\mapsto\mu^{\sig(x)}\phi$.
Conversely, suppose we know  the $\La^x$. To produce  $\apaircd$,
we only need to compute $\apair{\phi}{\psi}$  for $\phi,\psi\in\RR^d$,
since $\RR^j$ is orthogonal to $\RR^k$ if $j\neq k$.  By adjointness again we find
\begin{equation}\label{eq:apair_La} 
\apair{\mu^{x_1}\cdots\mu^{x_d}}{\psi}=
\La^{\sig(x_1)}\cdots\La^{\sig(x_d)}(\psi),
\qquad  \mbox{ if }\psi\in\RR^d
\end{equation}

The  cleanest formula for  $\apair{\phi}{\psi}$ comes from (\ref{eq:apair=}).
Let  $\bbT:\RR\to\C$ be the  projection operator  defined by the grading of $\RR$. 
Notice that $\bbT$ is \emph{classical},  i.e., we
know it before we quantize anything.
Recall the map $\phi\mapsto\phi\ssig$ from the proof 
of Proposition \ref{prop:RjRk}; this is also classical.  
$\bbT$ and $\T$ correspond  via $\bq$ and so  $\bbT$ is a $\circ$-trace 
by   Proposition \ref{prop:T}; 
we view this as the  ``quantum analog" of the fact that
$\bbT$  vanishes on Poisson brackets.   So (\ref{eq:apair=}) gives  
\begin{equation}\label{eq:Tcirc} 
\apair{\phi}{\psi}=\T\left(\bq(\phi)\bq(\psi\ssig)\right)=
\bbT(\phi\circ\psi\ssig),\qquad 
\phi,\psi\in\RR 
\end{equation}
For  $\phi,\psi\in\RR^d$, this reduces to
$\apair{\phi}{\psi}=\CR_d(\phi,\psi\ssig)$ where $\CR_p$ are the coefficients of
$\star$.

We can now characterize $\star$ without the explicit use of symmetry and unitarity.
\begin{prop}\label{prop:sat} 
The preferred star product $\star$ on $\RR$ we found in Theorem
\textup{\ref{thm:main}} is  uniquely determined by just the two properties:
\textup{(i)}  $\star$ corresponds to a $G$-equivariant quantization map
$\bq:\RR\to\cD$, and \textup{(ii)} $\star$ satisfies \textup{(\ref{eq:mux*phi})}
where the $\La^x$ are any operators.
\end{prop}
\begin{proof}
Suppose $\star$ satisfies  (i) and (ii). Then  $\star$ satisfies (\ref{eq:RjRk}) 
and so  $\bbT(\RR^j\circ\RR^k)=0$   for  $j\neq k$. 
Equivalently,    $\T(\cV^j\cV^k)=0$  if $j\neq k$.
We claim that this uniquely determines $\oplus_{d=0}^\infty\cV^d$
among all  $\g$-stable splittings  of   the order filtration of $\DD$. 
For it implies that  the spaces $\cV^d$ are orthogonal with respect to the  
symmetric   bilinear pairing  $\la(A,B)=\T(AB)$. But we know  $\la$ is 
non-degenerate on $\cV^d$; this follows because  $\cV^d$ is $\sig$-stable
and $\la(A,A\ssig)=\ga(A,A)$ is positive if $A\neq 0$. So there is only one
$\la$-orthogonal splitting. This proves our claim.
\end{proof}
We remark  that (\ref{eq:mux*phi}) and (\ref{eq:RjRk}) are in fact equivalent.

\begin{remark}\label{rem:theta}  
The  involution $\psi\mapsto\psi\ssig$ differs from pointwise complex conjugation
$\psi\mapsto\opsi$ precisely because $\GR$  is not  compact. But
there is a $\GR$-invariant \emph{indefinite} hermitian pairing on $\DD$ which may 
be more natural for  our quantization problem. This is   
$\tau(A,B)=\T(A\oB)$. The positivity of $\ga$ easily implies that there is a 
unique $\tau$-orthogonal splitting of the order filtration of $\DD$, which is again
$\oplus_{d=0}^\infty\cV^d$ .
So $\tau$ produces the same quantization map $\bq$. Then $\tau$ corresponds to 
the pairing $\pair{\phi}{\psi}=\bbT(\phi\circ\opsi)$. 
\end{remark}

\section{$\hatRR$ is a Fock space type model of $L^2(X,\Ehalf)$}  
\label{sec_fock}   
 
Combining the discussion in \S\ref{sec_thm} with our work in \S\ref{sec_I}, we find
\begin{cor}\label{cor:hatRR} 
The Hilbert space   completion $\hatRR=\widehat{\oplus}_{d=0}^\infty\RR^d$ 
of $\RR$ with respect to $\apaircd$ becomes
a   holomorphic  model for the unitary representation of $G$ on $L^2(X,\Ehalf)$.
We have, for the Harish-Chandra modules,  the  explicit intertwining isomorphism 
\begin{equation}\label{eq:long} 
\RR\mmapright{\bq}\DD\mmapright{\Delta}\GEhalffin
\end{equation}
\end{cor} 

While  ${L^2}(X,\Ehalf_X)$ is itself a Schroedinger type model, our
$\hatRR$ is a generalization of the Fock space  model of the oscillator \rep\ 
of the metaplectic group. This follows for three reasons. First, $\hatRR$ is the
completion of a space of ``polynomial" holomorphic functions. (We conjecture that
$\hatRR$  is a Hilbert space of holomorphic functions on $\TX$. This is true
when $\GR=SL_{n+1}(\R)$ and $\XR=\RP^n$ -- see \cite{A-B:gq}.)   

Second, the action of the  skew-hermitian operators $\pi^{x,\sig(x)}$
corresponding to the non-compact
part of $\gs$ is given by creation and annihilation operators.
For the non-compact part of $\gs$ is   $\{(ix,-ix)\,|\, x\in\gc\}$ and 
(\ref{eq:pixyphi=})  gives 
\begin{equation}\label{eq:noncompact} 
\pi^{ix,-ix}=2\mu^{ix}+2\La^{ix}
\end{equation} 
The   multiplication operators $\mu^{ix}$ are ``creation" operators mapping
$\RR^d$ to $\RR^{d+1}$, while the $\La^{ix}$  are ``annihilation" operators
mapping $\RR^d$ to $\RR^{d-1}$.  

Third,  the operators $\pi^{x,x}$ 
corresponding to the compact part $\{(x,x)\,|\, x\in\gc\}$ of $\gs$ 
are just the derivations   $\{\mu^{x},\cdot\}$ and these map $\RR^d$ to $\RR^d$.
Notice that the  operators $\mu^{ix}$ and  
$\{\mu^{x},\cdot\}$ are classical objects, while the   $\La^{ix}$ are  
quantum objects  (which encode $\apaircd$).

This gives new examples in the orbit method. For   $\RR$ identifies with the
algebra  of $G$-finite holomorphic functions on the
complex nilpotent orbit $\OO$ associated to $\TX$  (cf. \S\ref{sec_RR}).
We may regard $\OO$ as a real coadjoint orbit of  $G$, where
$\OO$ is enjoying a complex polarization. Then Theorem \ref{thm:main} and
Corollary \ref{cor:hatRR} give  a  quantization of  $\OO$.

\bibliographystyle{plain}

\end{document}